\documentclass[12pt]{amsart}
\usepackage[centertags]{amsmath}
\usepackage{euscript, mathrsfs}
\usepackage{amsmath}
\usepackage{amsfonts}
\usepackage{amscd}
\usepackage{amssymb}

\newcommand{\om}{{\omega}}

\newcommand{\Si}{{\Sigma}}

\newcommand{\ra}{\rightarrow }
\newcommand{\x}{\times}

\newtheorem{theorem}{Theorem}[section]
\newtheorem{thm}[theorem]{Theorem}
\newtheorem{cor}[theorem]{Corollary}

\newtheorem{rem}[theorem]{Remark}

\newtheorem{prop}[theorem]{Proposition}
\newtheorem{ex}[theorem]{Example}

\numberwithin{equation}{section}

\newtheorem{pr}[theorem]{Problem}
\newtheorem{conj}[theorem]{Conjecture}

\begin{document}
\title[Exotic smooth structures and symplectic forms]{Exotic smooth structures and symplectic forms on closed manifolds}
\author{Bogus\l\/aw Hajduk}
\author{Aleksy Tralle}
\address{Institute of Mathematics, Wroc\l\/aw University,
Plac Grunwaldzki 2/4, 50-384 Wroc\l\/aw, Poland}
\email{hajduk@math.uni.wroc.pl}
\urladdr{www.math.uni.wroc.pl/{~}hajduk}
\address{Department of Mathematics and Computer Science, University of Warmia and Mazury,
\.Zo\l\/nierska 14, 10-561 Olsztyn, Poland}
\email{tralle@matman.uwm.edu.pl} \urladdr{wmii.uwm.edu.pl/{~}tralle}


\begin{abstract}
In this paper we discuss relations between symplectic forms and
smooth structures on closed manifolds. Our main motivation  is the
problem if there exist symplectic structures on exotic tori. This is
a symplectic generalization of a problem posed by Benson and Gordon.
We give a short proof of the (known) positive answer to the original
question of Benson and Gordon that there are no K\"ahler structures
on exotic tori. We survey also  other related results which give
an evidence for the conjecture that there are no symplectic
structures on exotic tori.
\end{abstract}
\maketitle

\section{Introduction}
One of the most fundamental problems in symplectic topology is the
problem of existence of symplectic structures on closed manifolds.
This problem can be understood as follows. Consider an almost
complex $2n$-manifold $(M,J)$ with the "cohomological candidate for
the symplectic form", i.e. a class $a\in H^2(M;\mathbb{R})$ such
that $a^n\not=0.$ Then we can ask whether $M$ carries a symplectic
structure $\omega$ such that $[\omega]=a$ and $\omega$ is compatible
with an almost complex structure homotopic to $J$? This fundamental
question seems to be very difficult to answer in full generality.
Moreover, the work of many authors \cite{CFG, Cv, FeGG, FeS, FeM,
IRTU, McD, RT, TO} gives a strong evidence  that there are no
homotopic properties specific for symplectic manifolds, except the
obvious one given above. Various observations of particular classes
of symplectic manifolds are manifested in the {\it Thurston
conjecture: for any graded commutative finite-dimensional algebra
$H=\oplus_{i=0}^{2n}H^i$ satisfying the Poincar\'e duality and
possessing an
 element $a\in H^2$ such that $a^n\not=0$, there exists a closed
 symplectic manifold $(M,\omega)$ such that $H^*(M;\mathbb{R})\cong H$}
 (see \cite{Th}). In the sequel we will call such $H$ a {\it cohomologically
 symplectic} or $c$-symplectic algebra. If $M$ is a closed (not necessarily
 symplectic) $2n$-manifold with  $a\in H^2(M;\mathbb{R})$ such that $a^n\neq 0,$
  it is called a {\it $c$-symplectic manifold}. It is clear that to
 gain some understanding of the Thurston conjecture one needs at least examples of
 closed $c$-symplectic but {\it non-symplectic} manifolds.
\vskip6pt
 In the four-dimensional case examples of almost complex cohomologically
 symplectic but non-symplectic manifolds are given by the Seiberg-Witten
 theory \cite{T, FS, WG}. For example, it has become a classical result
 now that manifolds $k\mathbb{C}P^2\#l\overline{\mathbb{C}P^2}, k > 1$
 carry no symplectic structures compatible with the orientation given by
 the complex structure, although some of them are almost complex and they
 have the cohomology type of a closed symplectic manifolds. This follows,
 since their Seiberg-Witten invariant (which is a smooth invariant) vanishes,
 which contradicts the symplecticness by the Taubes theorem. There are several
 other results in this direction obtained by methods of the Seiberg-Witten
 theory (see, for example, \cite{FS,WG}). On the other hand, there are no
 known results of this type in higher
dimensions.

In view of \cite{Cv, FeM, IRTU, FeS, TO} on homotopy of symplectic
manifolds and of the crucial developments on smooth invariants in
the four-dimensional case one can ask how the existence of
symplectic structure does depend on smooth structure. Hence, the
following problem seems to be  very interesting. \vskip6pt

\begin{pr} \label{pr1} {\it It is known that for any $m\geq 5$, there are
smooth manifolds $\mathcal T^m$  which are homeomorphic but not
diffeomorphic to the standard torus $\mathbb{T}^m$. Does an exotic
torus $\mathcal T^{2n}, n>2,$ carry symplectic structures?}
\end{pr}
This is a generalization of a similar problem, posed by Benson and
Gordon in \cite{BG} for K\"ahler manifolds . In the present paper we
will give a short proof of the fact {\it that there are no K\"ahler
structures on exotic tori}. When this paper was written we have
learned that this fact was known to algebraic geometers \cite{BC,
C1}.

\vskip6pt In fact,  Baues, Cort\'ez and Catanese in these papers
developed a theory of aspherical K\"ahler manifolds with solvable
fundamental group. This theory has an application to another problem
of Benson and Gordon \cite{BG1}. A {\it solvmanifold} is a
homogeneous space $G/H$ with solvable Lie group $G$. In the sequel we always
assume that $H$ is a discrete cocompact subgroup (denoted by
$\Gamma$). We use the notation $L(G)$ for the Lie algebra of $G$. A
solvmanifold $G/\Gamma$ is called {\it completely solvable}, if for
any $x\in L(G)$ the linear operator $ad\,x: L(G)\rightarrow L(G)$
has only real eigenvalues.
\begin{pr}\label{bg2} Is every completely solvable solvmanifold $G/\Gamma$ which
carries a K\"ahler structure diffeomorphic to a torus?
\end{pr}
 As a result of efforts of \cite{A,AN,Ha,Ha1,TK} together with developments from
 algebraic geometry (see \cite{BC,C1}) the final (affirmative) solution of this
 problem was recently achieved (See Section 3).
 To the authors knowledge, the final clean proof was given in
 \cite{Ha, Ha1}. In \cite{C,C1, ABCKT} some properties of the Albanese map
 of K\"ahler $K(\pi, 1)$-manifolds with solvable $\pi$ were established. These
 properties allow one to show that no K\"ahler structures can exist on exotic
 tori, hence to answer the Benson-Gordon question in the negative.

Since we are interested in the symplectic version of the
 problem,  our approach goes in other direction, where explicit constructions of
 exotic differential structures play a role. Thus this paper can be
 considered as complementary to
\cite{BC, C1}, with the intersection  exactly in the original
Benson-Gordon problem.

   In \cite{HT} we obtained a result which gives a partial negative answer to
the symplectic question.  Sections 5 and 6 contain an exposition of it together 
with other related problems. 
This  provides some evidence for
the following conjecture.
\begin{conj}\label{conj}
There are no symplectic structures on exotic tori.
\end{conj}
The purpose of this article is to survey  results around this
conjecture as well as to discuss some relations between smooth
structures and symplectic forms on closed manifolds.

\section{Cohomologically symplectic manifolds}

A lot of work has been done by many mathematicians with the aim of better
understanding the homotopic properties of closed symplectic manifolds. These
results were initiated by Thurston's discovery that non-K\"ahlerness of
symplectic structures can be detected by homotopic invariants
\cite{BT, CFG, FeGG, FeS, G1, McD,RT, Th, TO}. It is understood now that
{\it all} known homotopic properties of K\"ahler manifolds may be violated by
closed symplectic manifolds. For example, it is known that K\"ahler manifolds
\begin{enumerate}
\item have even odd-dimensional Betti numbers;
\item satisfy the hard Lefschetz property;
\item have vanishing Massey products;
\item are formal.

\end{enumerate}
For a thorough discussion of the  homotopic properties (1)-(4) and
their role in K\"ahler and symplectic theory  we refer to the
monograph \cite{TO}. Examples of symplectic manifolds violating
(1)-(4) are constructed, for instance in \cite{Th, G1, McD, FeM, FeS,
BT, RT}. Moreover, it was asked what are the relations between
(1)-(3) in the following sense:
\begin{pr}\label{irtu} Can any combination of properties (1)-(3) or
their negatives be realized by a
closed symplectic manifold?
\end{pr}
A rather detailed answer is contained  in \cite{IRTU,FeM, Cv}. Results from these papers are summarized in the following statement.
\begin{thm}\label{ibeal} For any combination of properties (1)-(3) with the
following exceptions:
\begin{enumerate}
\item Massey products are trivial, the hard Lefschetz property holds and some
odd-dimensional Betti number is odd;
\item Massey products are non-trivial, the hard Lefschetz property holds and some
odd-dimensional Betti number is odd,
\end{enumerate}
there exists a closed symplectic manifold $(M,\omega)$ possessing this combination
of properties.
\end{thm}

Note that (1) and (2) are mentioned only for completness, since the
hard Lefschetz property implies the evennes of the odd-dimensional
Betti numbers.

The described results provide an  evidence for
 the Thurston conjecture. On the other hand there do exist relations between symplectic
structure and smooth structure on a closed manifold. This relation
is indicated by results of Taubes \cite{T}.
\begin{thm}\label{tb} Let $X$ be a compact, oriented, 4-dimensional manifold with
$b^+_2\geq 2$. Let $\omega$ be a symplectic form on $X$ with
$\omega\wedge\omega$ giving the orientation. Then the first Chern
class of the associated almost complex structure on $X$ has
the Seiberg-Witten invariant equal to $\pm 1$.
\end{thm}
This implies that connected sums of 4-manifolds with non
negative-definite intersection forms do not admit symplectic forms
which are compatible with the given orientation. The latter fact
follows, since Taubes also proves that if $X$ has $b^+_2\geq 2$ and
can be split by an embedded 3-sphere into $X_1\# X_2$ where neither
$X_1$ nor $X_2$ have negative definite intersection forms,
then the
Seiberg-Witten invariants of $X$ vanish. In particular, we have the
following corollary.
\begin{cor}\label{cp2} Connected sums $k\mathbb{C}P^2\# l\overline{\mathbb{C}P}^2$
when $k>1$ do not admit symplectic forms  which define the
orientation given by the almost complex structure.
\end{cor}
The Ehresmann-Wu theorem which characterizes 4-manifolds  admitting
almost complex structures \cite{Be} easily implies that
$k\mathbb{C}P^2\# l\overline{\mathbb{C}P}^2$ is almost complex if
and only if  $k$ is odd.
\vskip6pt
Another kind of examples can be obtained by the Seiberg-Witten theory as follows.
\begin{thm}\label{wg}{\rm (\cite{WG})} Let $\tilde X$ be a K\"ahler 4-manifold
with positive basic class $K_{\tilde X}>0$ and $b_2^+(\tilde X)>3$. Suppose that
$\sigma: \tilde X\rightarrow\tilde X$ is an antihilomorphic involution without
fixed points. Then for $X=\tilde X/\sigma$, $SW(X)=0$, and, hence, $X$ cannot
carry symplectic structures.
\end{thm}
Here is an explicit example.
\begin{ex}\label{e-wg} {\rm Let $\tilde X$ denote the hypersurface
$$\tilde X=\{[z_0:z_1:z_2:z_3]\in\mathbb{C}P^3\,
\,|\,\sum_{i=0}^3z_i^{2d}=0,\,d>2\}\subset\mathbb{C}P^3.$$
Let $\mathbb{Z}^2$ act on $\tilde X$ by complex conjugation.
By Theorem 2.3\label{wg}, $X=\tilde X/\mathbb{Z}^2$ has no symplectic structure.}
\end{ex}
Seiberg-Witten invariants yield also more sophisticated examples,
for instance, in the simply-connected case. This is the {\it knot
construction} due to Fintushel and Stern (see \cite{FS} and
\cite{GS}). The relation between symplectic structures and
Seiberg-Witten invariants is delicate, since there exist closed
almost complex manifolds with Seiberg-Witten invariants satisfying
the Taubes condition and having no symplectic structures \cite{KMT}.
\vskip6pt

In contrast to dimension 4, the following question is open.
\begin{pr}\label{pr2} Do there exist closed $c$-symplectic non-symplectic (almost
complex) manifolds in dimensions $\geq 6$?
\end{pr}
In fact, up to now, no general tools have been worked out. A rich algebraic theory
being a source of potential examples was created by Lupton and Oprea in \cite{LO}.
However, although the  examples available from \cite{LO} are $c$-symplectic,  there is no
 way in sight neither to prove the existence of symplectic structures nor to disprove this.
To give the reader some taste of what is known in higher dimensions we recall here an
example from \cite{LO}.
\begin{ex}\label{lo} {\rm Let $K=\mathbb{C}P^2\times V$, where $V$ is a hypersurface
in $\mathbb{C}P^4$ defined by a single equation of degree 3. The rational cohomology
algebra of $K$ can be easily calculated, since the cohomology of hypersurfaces in
$\mathbb{C}P^n$ is known. In the case considered here we have
$$H^*(K;\mathbb{Q})\cong \mathbb{Q}[\omega]/(\omega^3)\otimes
\mathbb{Q}[x,a_1,...,a_5,a^*_1,...,a_5^*]/\mathcal{R},$$
with the ideal $\mathcal{R}$ generated by
$$\{xa_j,xa_j^*\}_{j=1,...,5},\,\{a_ja_k,a^*_ja^*_k\}_{1\leq j<k\leq 5},
\{a_ja^*_k\}_{j\not=k}$$
 and
$$x^3-(\sum_{j=1}^5a_ja^*_j).$$

It is shown that there exists a {\it non-formal} space $X$ such that
$H^*(X)\cong H^*(K)$. Note that $K$ is K\"ahler, and it is  well
known that K\"ahler manifolds are formal. Using rational surgery
\cite{Ba}, Lupton and Oprea show that there exists a simply
connected 10-dimensional manifold with the same minimal model, and
cohomology as $X$ (and $K$) (see \cite{TO} for the terminology).
Finally, we have arrived at the conclusion:} {\it there exists a $c$-symplectic simply connected 
10-dimensional smooth
manifold which  has the cohomology algebra of a K\"ahler manifold and homotopy type different from that of K\"ahler
manifolds}.
\end{ex}
\section{Albanese map  and a solution of the Benson-Gordon problem}
The theorem below yields a solution to the original problem \cite{BG}.
Although it is an evidence to our conjecture, since K\"ahler forms are
symplectic,  the proof is "non-symplectic" and probably cannot be generalized.
\begin{thm}\label{th1} There are no K\"ahler structures on exotic tori.
\end{thm}
\vskip6pt
The proof of this result uses some properties of the Albanese
map, which we recall now.
Let $X$ be a compact K\"ahler manifold.  By definition
(see \cite{BPV}), the {\it Albanese variety} of $X$ is the complex torus
$$Alb(X)=H^0(X,\Omega_X^1)^*/j(H_1(X,\mathbb{Z})),$$
where $j$ is the homomorphism
$$j: H_1(X,\mathbb{Z})\to H^0(X,\Omega_X^1)^*,$$
$$j([\gamma])=(\omega\to\int_{\gamma}\omega).$$
Here $\Omega_X^1$ denotes the sheaf of germs of holomorphic 1-forms
on $X$. It is known (see \cite{ABCKT}, \cite{BPV}), that
$j(H_1(X,\mathbb{Z}))=H_1(X,\mathbb{Z})/Torsion$ is a lattice of
rank $b_1(X)=2h^0(X,\Omega_X^1)$, and thus $Alb(X)$ is a complex
torus. Fixing a basepoint $x_0\in X$, one defines the Albanese map
by
$$\alpha_X: X\to Alb(X),\,\alpha_X(x)=(\omega\to \int_{x_0}^x\omega).$$
The following properties of the Albanese map are known (\cite{ABCKT}, \cite{BPV}):
\begin{itemize}
\item $\alpha_X$ is a holomorphic map from $X$ to the complex torus
$\mathbb{T}^{b_1(X)}=Alb(X)$, of complex dimension
${\frac12}b_1(X)$;
\item
  $\alpha_X$ induces a surjection $(\alpha_X)_*:\pi_1(X)\to \pi_1(Alb(X))$;
\end{itemize}
For any topological space $X$ of finite type define $a(X)$
to be the maximal integer for which the image of
$\Lambda^mH^1(X,\mathbb{R})$ in $H^m(X,\mathbb{R})$ is non-trivial.
The following result can be found in \cite{ABCKT, C}.
\begin{thm}  Let $X$ be a compact K\"ahler manifold. Then
$a(X)$ is the real dimension of its Albanese image:
$$a(X)=\dim\alpha_X(X).$$
\end{thm}
Having the above facts in mind we are ready to prove Theorem 3.1.

\vskip6pt {\it Proof of Theorem 3.1.} Let $X$ denote a smooth
manifold homeomorphic to $\mathbb{T}^{2n},$ the torus of dimension
$2n$. Clearly, $a(X)=2n,$ since $\Lambda^{2n}H^1(X,\mathbb{R}) \cong
H^{2n}(X,\mathbb R)$ is spanned by $x_1\wedge ...\wedge
x_{2n}\not=0$, where $x_i$ denote the 1-dimensional generators of
$H^1(X,\mathbb{R})$.

\vskip6pt Assume that $X$ carries a K\"ahler structure. Then we have
the Albanese map $\alpha_X:X\to \mathbb T^{2n}.$

Consider a regular point $p$ in the image of $\alpha_X(X).$ Since
holomorphic map preserves orientations, its degree is equal to the
cardinality of $\alpha_X^{-1}(p).$ As $\dim\alpha_X(X)=2n=\dim
Alb(X),$ it implies that $\deg\alpha_X>0$ and $\alpha_X$ is onto.
Moreover, $\alpha_X$ induces a surjection of free abelian groups of
rank $2n,$
$$(\alpha_X)_*:\pi_1(X)\to \pi_1(\mathbb{T}^{2n}),$$
hence an isomorphism. This gives also that $\alpha_X$ induces an
isomorphism  $H^1(\mathbb T^{2n};\mathbb Z)\to H^1(X;\mathbb Z)$ and
consequently an isomorphism $$H^{2n}(\mathbb T^{2n};\mathbb Z)\to
H^{2n}(X;\mathbb Z).$$ Thus $\alpha_X$ is a map of degree one.

Since the set of regular points of a holomorphic onto map is dense,
we have that $\alpha_X$ is bijective on a dense subset, thus everywhere.
It is a standard fact that a holomorphic homeomorphism is
biholomorphic (cf. \cite{FG}, Th. 8.5), thus $\alpha$ is
a diffeomorphism.
The proof is complete.
\begin{rem}\label{bau-c} {\rm As we have mentioned in the introduction,
Baues and Cort\'ez \cite{BC} gave an account on recent developments concerning
the classification of compact aspherical K\"ahler manifolds whose fundamental
groups contain a solvable subgroup of finite index. These developments lead
to solutions of the both Benson-Gordon problems. In fact, more general results hold.}
\begin{thm}\label{t1bc}
Let $M$ be an infra-solvmanifold which admits a K\"ahler metric. Then $M$ is
diffeomorphic to a flat Riemannian manifold.
\end{thm}
\begin{thm}\label{t2bc} Let $X$ be a compact K\"ahler manifold of complex
dimension $n$. Assume that $X$ satisfies $2\dim\,X=\dim\,H^1(X;\mathbb{C})$.
If $\dim\,H^2(X;\mathbb{C})\leq {2n\choose 2},$ and $H^{2n}(X;\mathbb{Z})$ is
generated by integral classes of degree 1, then the Albanese morphism
$\alpha_X: X\to \operatorname{Alb}(X)$ is a biholomorphic map.
\end{thm}
{\rm The proof of the last theorem is the same as Theorem
3.1\label{tm1}, which follows then as a corollary. The same result
follows from \cite{C1}. Proposition 4.8 in that paper implies
Theorem 3.1, with the proof going along the same lines, although the
solution of the Benson and Gordon problem is apparently a byproduct,
and is not mentioned.}

\end{rem}
\section{Examples of exotic structures on tori}
Various examples of different smooth structures within the homotopy type of
manifolds can be found in \cite{HS,HW,W}. However, our aim is to present
explicit constructions which enable to formulate our problems in the language
of groups $\operatorname{Diff}\,(M)$ and $\operatorname{Symp}\,(M,\omega)$.
The simplest examples of exotic tori are obtained as connected
sums of the standard tori with  homotopy spheres. Namely, consider a
homotopy sphere $\Si$ of dimension $k$ and the manifold $\mathcal
T_{\Si} = \mathbb T^k \# \Si \x \mathbb T^{2n-k}.$

 The
connected sum operation can be understood as follows. Let $\Si= D^k\cup_fD^k,$ when $f \in \operatorname{Diff}(S^{k-1},S^{k-1}_-),$ which means
that we glue two copies of the disk using a diffeomorphism of the
boundary sphere supported in the upper half-sphere. Then $\mathbb T^k \#
\Si$ is obtained by cutting $\mathbb T^k$ along an embedded
$(k-1)$-disk and gluing it again along the disk using $f.$ Up to
an orientation choice (which replaces $\Si$ by $-\Si ),$ this does not
depend on the choice of the disk, thus we can choose it in a
subtorus $\mathbb T^{k-1}\subset \mathbb T^k.$ This is equivalent to
cutting $\mathbb T^k$ along the subtorus and gluing again with a
diffeomorphism $\hat f \in \operatorname{Diff}(\mathbb T^{k-1})$ extending $f$ by
the identity. Thus what we get the mapping torus of $\hat f,$ i.e., the fibration over the circle with
fiber $\mathbb T^{k-1}$ and gluing diffeomorphism $\hat f,$ (cf.
\cite{H}).

 For $\mathcal T = (\mathbb T^k \# \Si )\x \mathbb
T^{2n-k}$ we get a fibration over $\mathbb T^{2n-k+1}$ and fiber
$\mathbb T^{k-1}.$ If $k$ is odd, there exist symplectic structures
on base and fiber, so we can ask if there exist such symplectic
structures which induce a symplectic structure on $\mathcal T.$
Since the fibration is homotopically (and even topologically)
trivial, a theorem of Thurston \cite{MS} says that this is the case if $\hat f$
is isotopic to a symplectomorphism (with respect to some symplectic
structure on the fiber). The main purpose of \cite{HT} was to show
that, in general, there is an obstruction to such isotopy.

In general, exotic differential structures on tori are obtained by
iterating the operation described above, but then we have to apply
the cutting and pasting along possibly exotic subtori. It follows
from \cite{W}, "Fake tori" chapter, that the examples above give in
fact all nonstandard differential structures on tori.
However, it is not easy to distinguish diffeomorphism type of two
manifolds obtained by the construction above. The classification
resulting from surgery or from smoothing theory identifies smooth
structures up to diffeomorphisms isotopic topologically to the
identity. Here we need to know when two structures are
diffeomorphic, thus we have to see what is the action of $\pi_0\operatorname{Diff}(\mathbb T^k)$
on the set of smoothings. Obviously, diffeomorphisms of the torus
act nontrivially on homology, hence regluing along two subtori of
the same dimension and using the same homotopy sphere yields
diffeomorphic structures. Simply one can exchange the subtori by a
 diffeomorphism and it gives a diffeomorphism of resulting
manifolds. However, $\pi_0\operatorname{Diff}(\mathbb T^k)$ is unknown for $k>4,$
which makes further  analysis difficult.

The simplest case, but still nontrivial, is that of the connected
sum. It follows from \cite{K} that $\mathbb T \# \Si_1 \cong
\mathbb T \# \Si_2$ if and only if $\Si_1 \cong \Si_2.$

If we need only some examples of non diffeomorphic exotic tori, one
can use the Atiyah - Milnor - Singer invariant $\hat a$ which
distinguishes cobordism classes of framed manifolds. In fact we need
only its $\mathbb Z_2$ part which is known as Hitchin invariant. The
instructive example is $k=2n-1=8s+1$. We have $\mathcal T = \mathbb
T^{8s+1} \# \Si\x \mathbb T^1.$

 In dimension $8s+1$ and $8s+2$ half of
homotopy spheres have nontrivial generalized $\hat a$-genus
\cite{H}. The $\hat a$-genus can be defined, for any closed spin
manifold $M^m$, as the $KO$-theoretical index of the Dirac operator
with values in $KO^{-m}(pt).$ It is known that the coefficient
groups $KO^{-*}(pt)$ are the following
$$
KO^{-m}(pt)=\begin{cases}
\mathbb{Z} &\text{for}\ \  m\equiv 0\ (\text{mod}\, \ 4);\\
\mathbb{Z}_2  &\text{if}\ \  m\equiv 1,2\ (\text{mod}\,\ 8);\\
0 &\text{for any other}\,\, m.
\end{cases}
$$
Let $f: M^m\to \{pt\}$ denote the obvious collapsing map. For a spin
structure on $M^m$ we have the Gysin map $f_{!}:KO^0(M^m)\to
KO^{-m}(pt).$ By definition, the {\it $\hat a$-genus} of $M^m$ is an
element of $KO^{-m}(pt)$ given by the formula
$$\hat a(M)=f_{!}(1).$$

The ${\hat a}$-genus has the following properties (see \cite{LM}):
\begin{enumerate}
\item for any closed spin manifolds $X$ and $Y$

$$ \hat a(X\#Y) = \hat a(X) + \hat a(Y),\ \
(\text{when}\, \dim\, X = \dim\, Y)$$ and $$ \hat a (X \x Y) = \hat
a(X) \hat a(Y), $$

\item $\hat a$ is a spin cobordism invariant,

\item for any $m > 2$ and any spin structure on the standard torus
$\mathbb T^m$ we have $\hat a(\mathbb T^m) = 0.$
\end{enumerate}
Note that (3) follows from (1) and (2) since any spin structure
on $\mathbb T^m$ is given as a product of spin structures on circles.
One of the two possible structures on $S^1$ has nonzero $\hat
a$-genus, but the third power of the nonzero element of
$KO^{-1}(pt)$ vanishes.

\begin{prop}
The manifold $\mathcal T = (\mathbb T^{8s+1}\# \Si)\x S^1$ is
homeomorphic, but not diffeomorphic to the standard torus $\mathbb
T^{8s+2}$ if $\hat a(\Si ) \neq 0.$ Moreover, the $\hat a$-genus of
$\mathcal T$ does depend on the choice of the spin structure.
\end{prop}

By properties (1) and (2), if we choose the trivial spin structure on
$\mathbb T^{8s+1}$ and nontrivial one on $S^1,$  one has $\hat
a(\mathcal T)\neq 0,$ while it is always zero for the standard
torus. This argument works also for $k=8k+1$ or $8k+2$ and $2n-k>1.$
One may consider the map $A: H^1(\mathcal T; \mathbb Z^{2n-k})\ra
\mathbb Z_2$ given as follows. Let $\phi : \mathcal T \ra \mathbb
T^{2n-k}$ correspond to $x\in H^1(\mathcal T; \mathbb Z^{2n-k})$
(note that $\mathbb T^r$ is the Eilenberg-MacLane space $K(\mathbb
Z^r, 1)).$ Let $A(x) = \hat a(\phi^{-1}(p)),$ where $p$ is a regular
value of $\phi .$  In our examples there exists a spin structure
such that $A$ is nontrivial, but it is trivial for the standard
torus and $k>2.$

Finally, it is known that all homotopy tori are parallelizable (see
\cite{W} or \cite{HT}). Thus any such manifold admits an almost
complex structure and it is  cohomologically symplectic. Therefore,
an exotic torus with no symplectic structure would be an example of
a non-symplectic manifold satisfying the two "obvious" conditions
necessary for symplecticness.

It is remarkable that our assumption on dimension plays a role in
the conjecture about the non-existence of symplectic structures. In
dimension 4 the situation is different \cite{P1, P2,P3}. Using the
Gompf symplectic sum construction \cite{G1, MW}, it is possible to
construct families of closed 4-manifolds which are homeomorphic to
$(2m+1)\mathbb{C}P^2\# n\overline{\mathbb{C}P}^2$ but {\it not
mutually diffeomorphic}. It follows that $(2m+1)\mathbb{C}P^2\#
n\overline{\mathbb{C}P}^2$ admits infinitely many symplectic
structures which {\it are all "exotic"}, since this manifold does
not admit symplectic structures coming from the standard smooth
structure, by the Taubes theorem. Let us give an exact formulation
of the main result in \cite{P1}.
\begin{thm}\label{park} For each positive integer $m$ and $n$ satisying
$2m+8\leq n\leq 10m+9$ there exists a family of simply connected, closed,
nonspin, irreducible symplectic 4-manifolds
$$\{X_{2m+1,n}(p)\,|\, p\,\text{is a positive integer}\,\}$$
which are all homemomorphic to $(2m+1)\mathbb{C}P^2\#
n\overline{\mathbb{C}P}^2$, but not mutually diffeomorphic.
\end{thm}

\section{Diffeomorphisms of tori and non-existence of fiberwise symplectic
structures on exotic tori}

In this section we will describe the results of \cite{HT}. We
consider the following question:

\begin{pr} Is there a symplectic structure on $\mathcal T =(\mathbb T^{8s+1}\# \Si)\x S^1$ compatible with 
the fiber bundle structure
$$\mathbb T^{8s}\rightarrow \mathcal T=(\mathbb T^{8s+1}\#\Sigma)\times S^1 \rightarrow \mathbb T^2?$$

\end{pr}

The necessary and sufficient conditions for a fibration over a
symplectic manifold with a symplectic fiber to have a symplectic
structure such that each fiber is symplectic was given by Thurston \cite{MS}:
one has to know that the fibration is symplectic and a cohomological
condition should be satisfied.  Clearly, the cohomology ring
$H^*(\mathcal T)$ is isomorphic to $H^*(\Bbb T^{2n})\otimes
H^*(\Bbb T^2)$. This easily implies that the cohomological condition
of Thurston's construction is  satisfied. The symplecticness of the
fibration is that the gluing diffeomorphism $f\in \operatorname{Diff}(\mathbb
T^{8s})$ must be isotopic to a symplectomorphism.

Thus we come to the following question.

\begin{pr} Given a diffeomorphism $f:\Bbb T^{8s}\to\Bbb
T^{8s}$ supported in  an embedded disc but non-isotopic to the
identity, is there a symplectomorphism in the isotopy class of $f$?
\end{pr}

We have seen above that the positive answer to Problem 5.2 would imply
that $\mathcal T$ admits a symplectic structure compatible with the
fibration. Our goal is to give negative examples to this problem in the particular case of the standard symplectic structure on $\mathbb T^{2n}$.   
 Let
$\pi_0(\operatorname{Diff_+}\,(M))$ denote the group of isotopy
classes of orientation preserving diffeomorphisms of a smooth
oriented manifold $M$. Assume now that $M$ is $2n$-dimensional and
admits almost complex structures, and let ${\Bbb J}M$ denote the set
of homotopy classes of such structures, compatible with the given
orientation. Any diffeomorphism $f$ acts on the set of all almost
complex structures by the rule
$$f_*J=df J df^{-1},$$
where $df: TM\to TM$ denotes the differential of $f$. This action
clearly descends to the action of
$\pi_0(\operatorname{Diff_+}\,(M))$ on ${\Bbb J}M$.

Let now ${\frak G}\, (M)$ denote the subgroup of
$\pi_0(\operatorname{Diff}\,(M))$ generated by diffeomorphisms with
supports in discs.

In the sequel we will show that there exist diffeomorphisms $f: \mathbb T^{8k}\rightarrow\mathbb T^{8s}$ supported in a disc, whose isotopy classes $[f]\in {\frak G}\,(M)$ do not preserve the homotopy class $[J_0]\in{\mathbb J}M$ of the standard complex structure. Therefore, they cannot be isotopic to symplectomorphisms with respect to the standard symplectic structure $\omega_0$. Indeed, any symplectomorphism carries any
almost complex structure compatible with a symplectic form to
another almost complex structure compatible with the same symplectic
form, but the space of all such almost complex structures is
contractible. 

Now, let us give a sketch of the main line of proof that such $f$ exist. We give
first a necessary homotopic condition on a diffeomorphism to be isotopic to a symplectomorphism.

\begin{thm} Let $f\in\operatorname{Diff}(\mathbb T^{4n})$ be supported in a disc $D^{4n}\subset\mathbb T^{4n}$. If $f$ is isotopic to a symplectomorphism with respect to the standard symplectic structure, then $df$ restricted to its support disc $D^{4n}$ gives in $\pi_{4n}SO(4n)$ the trivial homotopy class.
\end{thm}
 Since the $\hat a$-genus is additive
with respect to the operation of connected sum and it is nontrivial
for some homotopy spheres in dimension $8k+1,$ one easily concludes
that there exist isotopy classes of diffeomorphisms $f$ with support
in a disc such that $\hat a(\mathbb T^{8k+1}_f)\ne 0$. Here $\mathbb
T^{8k+1}_f$ denotes the mapping torus of $f$ (see Section 4).  On the other hand, we
prove the following result.

\begin{thm} In the notation of Theorem 5.3, if $[df]=0$, then
$\hat a(\mathbb T_f)=0$.
\end{thm}
Comparing Theorem 5.3 and Theorem 5.4 we get the conclusion.

\begin{thm} For any $k>0$ there exist diffeomorphisms
$f:\Bbb T^{8k}\to \Bbb T^{8k}$ with support in a disc which are not isotopic to a symplectomorphism of $(\mathbb T^{8k},\omega_0)$.
\end{thm}

\section{Further questions and perspectives}

We have shown in the preceding section that our question on
existence of symplectic structure on exotic tori is  related to the
problem  whether a diffeomorphism of an even dimensional torus,
supported in a disk and non-isotopic to the identity, can be
isotopic to a symplectomorphism.

Some more questions of that type were considered in symplectic
topology and give additional evidence for Conjecture 1.2. A similar
rigidity question was posed by McDuff and Salamon \cite{MS}: {\it Is
it true that any symplectomorphism of a torus which induces identity on
homology is isotopic to the identity?}

In \cite{CLO} the action of $\pi_4SO(4)$ on almost complex
structures on 4-manifolds was considered. The action is defined by
changing an almost complex structure in a disk by a map
$(D^4,\partial D^4)\rightarrow (SO(4),id)$ representing an element
of $\pi_4SO(4).$ The main result is that if an almost complex
structure is compatible with a symplectic structure, than the almost
complex structure  twisted by a nontrivial element such that the
corresponding $spin^c-$structure changes, has no compatible
symplectic structure.

Another related question is the following. Consider the space of
symplectomorphisms of $\mathbb R^{2n}$ with compact supports. It was
proven by Gromov \cite{G} that for $n=2$ this space is contractible.
The higher dimensional case is open. From what we have said before, a
weaker question seems to be natural: {\it is any symplectomorphism with
compact support in $\mathbb R^{2n}$ isotopic (smoothly) to the
identity?} Or, more generally, {\it is the image of the space of compactly
supported symplectomorphisms in the group of all compactly supported
diffeomorphism contractible?}

The latter question (at least for $\pi_0$) can be attacked along the
following lines motivated by \cite{MS2}. Let $f$ be such a symplectomorphism and $\om$ be
the standard symplectic structure. Then $f^*\om =\om,$ and for an
almost complex structure $J$ compatible with $\om$ we have $f^*J$
also compatible with $\om .$ Thus there exists a path $J_t$ of
almost complex structures compatible with $\om$ connecting $J$ with
$f^*\om .$ Because of the compact support assumption we can consider
all this in $\mathbb T^{2n}$ or in $\mathbb
T^{2n-2}\x S^2.$ Then the path of the evaluation maps for the space
of $J_t$-holomorphic curves should give an isotopy between $f$ and
the identity. We plan to address this problem in a forthcoming
paper.

 {\bf Acknowledgment.} We thank Jarek K\c edra for valuable
discussions. This work was supported, partially, by the Polish
Ministry of Science and Higher Education, research project no. 1P03A
03330, and by the Mathematisches Forschungsinstitut Oberwolfach
under the "Research in Pairs" program. The authors cordially thank
this institution for hospitality and excellent working conditions.


\begin{thebibliography}{ABCDE}
\bibitem[A]{A} D. Arapura, {\it K\"ahler solvmanifolds}, Internat. Math. Res. Notices 3(2004), 131-137.
\bibitem[AN]{AN} D. Arapura, M. Nori {\it Solvable fundamental groups of algebraic varieties and K\"ahler manifolds}, Compositio Math. 116(1999), 173-188.
\bibitem[ABCKT]{ABCKT} J. Amoros, F. Bogomolov, K. Corlette,
D. Kotschik, D. Toledo, {\bf Fundamental Groups of Compact K\"ahler
Manifolds}, Amer. Math. Soc., Providence, 1996.
\bibitem[Ba]{Ba} J. Barge, {\it Structures differentiables sur les types d'homotopie rationelle simplement connexes}, Ann. Sci. \'Ecole Norm. Sup. 9(1976),469-501.
\bibitem[Be]{Be} A. Besse, {\bf G\'eometrie riemannienne en dimension 4}, Cedic, Paris, 1981.
\bibitem[BC]{BC} O. Baues, V. Cort\'ez, {\it Aspherical K\"ahler manifolds with solvable fundamental group}, math.DG/0601616
\bibitem[BG]{BG}  C. Benson and C. Gordon, {\it K\"ahler
and symplectic structures on nilmanifolds}, Topology, {\bf 27}
(1988), 513-518.
\bibitem[BG1]{BG1} C. Benson, C. Gordon, {\it K\"ahler structures on compact solvmanifolds}, Proc. Amer. Math. Soc. 108(1990), 971-980.
\bibitem[BT]{BT} I. Babenko, I. Taimanov, {\it On nonformal simply connected manifolds}, Siberian Math. J. 41(2000),204-217.
\bibitem[BPV]{BPV} W. Barth, C. Peters, A. Van de Ven,
{\bf Compact Complex Surfaces}, Springer, Berlin, 1984.
\bibitem[C]{C} F. Catanese, {\it Moduli and classification of
irregular K\"ahler manifolds (and algebraic varieties) with Albanese
general type fibrations}, Invent. Math. {\bf 104} (1991), 263-289.
\bibitem[C1]{C1} F. Catanese, {\it Deformation types of real and complex manifolds},
math.AG/0111245
\bibitem[Cv]{Cv} G. R. Cavalcanti {\it The Lefschetz property, formality and blowing
up in symplectic geometry}, Trans. Amer. Math. Soc. 359(2007), 338-348.
\bibitem[CFG]{CFG} L.A. Cordero, M. Fern\'andez, A. Gray, {\it Symplectic manifolds
with no K\"ahler structure}, Topology 25(1986), 375-380
\bibitem[CLO]{CLO} F. Connolly, Le Hong Van, K. Ono,{\it Almost complex structures
which are compatible  with K\"ahler structures or symplectic structures}, Annals
Global Anal. Geom. 15(1997), 325-334.
\bibitem[FeGG]{FeGG} M. Fern\'andez, M. Gotay, A. Gray, {\it Compact parallelizable
four dimensional symplectic and compact solvmanifolds}, Proc. Amer.
Math. Soc. 103(1988), 1209-1212.
\bibitem[FeM]{FeM} M. Fern\'andez, V. Mu\~noz, {\it Formality of Donaldson submanifolds},
Math. Z. 250(2005), 149-175.
\bibitem[FeS]{FeS} M. Fern\'andez, M. de Le\'on, M. Saralegui, {\it A six-dimensional
compact symplectic solvmanifold without K\"ahler structures}, Osaka J. Math. 33(1996), 19-35.
\bibitem[FS]{FS} R. Fintushel, R. Stern, {\it Knots, links and 4-manifolds},
Invent. Math. 134(1998), 363-400.
\bibitem[FG]{FG} K. Fritzsche, H. Grauert, {\bf From Holomorphic
Functions to Complex Manifolds}, Springer, Berlin, 2002.
\bibitem[G]{G} M. Gromov, {\it Pseudoholomorphic curves in symplectic
manifolds,} Invent. Math. 82 (1985), 307 - 345.
\bibitem[G1]{G1} R. Gompf, {\it A new construction of symplectic manifolds},
Annals Math. 142(1995), 527-595.
\bibitem[GS]{GS} R. Gompf, A. Stipsicz, {\bf 4-Manifolds and Kirby Calculus},
Amer. Math. Soc., 1999.
\bibitem[H]{H} N. Hitchin, {\it Harmonic spinors}, Adv. Math. {\bf 14}(1974),
1-55.
\bibitem[Ha]{Ha} K. Hasegawa, {\it A note on compact solvmanifolds with
K\"ahler structures}, Osaka J. Math.43(2006), 131-135.
\bibitem[Ha1]{Ha1} K. Hasegawa, {\it Complex and K\"ahler structures on compact
solvmanifolds}, J. Symplectic Geom. 3(2005), 749-767.
\bibitem[HS]{HS} W.C. Hsiang, J.L. Shaneson, {\it Fake tori, the annulus
conjecture, and the conjectures of Kirby}, Proc. N.A.S. 62(1969), 687-691.
\bibitem[HW]{HW} W.C. Hsiang, C.T.C. Wall, {\it On homotopy tori II}, Bull.
London Math. Soc. 1(1969), 341-342.
\bibitem[HT]{HT} B. Hajduk, A. Tralle, {\it Diffeomorphisms and
almost complex structures on tori}, Annals Global Anal. Geom. {\bf
28} (2005), 337-349.
\bibitem[IRTU]{IRTU} R. Ib\'a\~nez, Y. Rudyak, A. Tralle, L. Ugarte,
{\it On certain geometric and homotopy properties of closed symplectic manifolds},
Topology Appl. 127(2003), 33-45
\bibitem[K]{K} A. Kosi\'nski, {\it On the inertia group of $\pi$-manifolds},
Amer. J. Math. 89(1967), 227-248.
\bibitem[KMT]{KMT} D. Kotschik, J. Morgan, C. Taubes, {\it Four-manifolds
without symplectic structures but with non-trivial Seiberg-Witten invariants},
Math. Res. Letters 2(1995),119-124.
\bibitem[LM]{LM} H. Lawson, M.L.Michelson, {\bf Spin Geometry}, Princeton
University Press, Princeton, 1989.
\bibitem[McD]{McD} D. McDuff, {\it Examples of simply-connected symplectic
non-K\"ahlerian manifolds}, J. Differential Geom. 20(1984),267-277.
\bibitem[MS]{MS} D. McDuff, D. Salamon, {\bf Introduction to
Symplectic Topology,} Oxford University Press, Oxford, 1998.
\bibitem[MS2]{MS2} D. McDuff, D. Salamon, {\bf $J$-holomorphic Curves and
Symplectic Topology}, Amer. Math. Soc., Providence RI, 2004.
\bibitem[MW]{MW} J. McCarthy, J. Wolfson, {\it Symplectic normal connect sum},
Topology 33(1994), 729-764.
\bibitem[LO]{LO} G. Lupton, J. Oprea, {\it Symplectic manifolds and formality}
J. Pure Appl. Algebra 91(1994),193-207.
\bibitem[P1]{P1} J. Park, {\it Exotic smooth structures on 4-manifolds},
Forum Math. 14(2002), 915-929.
\bibitem[P2]{P2} J. Park, {\it Exotic smooth structures on
$3\mathbb{C}P^2\# n\overline{\mathbb{C}P}^2$}, Proc. Amer. Math. Soc. 128(2000), 3057-3065.
\bibitem[P3]{P3} J. Park, {\it Exotic smooth structures on
$3\mathbb{C}P^2\# n\overline{\mathbb{C}P}^2$, Part II}, Proc. Amer. Math.
Soc. 128(2000), 3067-3073.
\bibitem[RT]{RT} Y. Rudyak, A. Tralle, {\it On Thom spaces, Massey products
and nonformal symplectic manifolds}, Internat. Math. Res. Notices, 10(2000), 495-513.
\bibitem[T]{T} C.H. Taubes, {\it The Seiberg-Witten invariants
and symplectic forms}, Math. Res. Letters, {\bf 1} (1994), 809-822.
\bibitem[Th]{Th} W.P. Thurston, {\it Some simple examples of symplectic manifolds},
Proc. Amer. Math. Soc. 55(1976), 467-468.
\bibitem[TK]{TK} A. Tralle, J. K\c edra, {\it Completely solvable K\"ahler
solvmanifolds are tori}, Internat. Math. Res. Notices 15(1997), 727-732.
\bibitem[TO]{TO} A. Tralle, J. Oprea, {\bf Symplectic Manifolds with no
K\"ahler Structure}, Lect. Notes Math. 1661, Springer, Berlin, 1997
\bibitem[W]{W} C.T.C. Wall, {\bf Surgery on Compact Manifolds},
2nd edn, Amer. Math. Soc., Providence, 1999.
\bibitem[WG]{WG} S. Wang, {\it A vanishing theorem for the
Seiberg-Witten invariants}, Math. Res. Letters, {\bf 2} (1995),
305-310.
\end{thebibliography}
\end{document}